\documentclass[letterpaper, 10 pt, conference]{ieeeconf}  

\IEEEoverridecommandlockouts                              
\usepackage{amsfonts}
\usepackage{epsfig}
\usepackage{subfigure}

\def\Dl{\Delta}

\def\lz{{l_0}}

\def\F{{\cal F}}
\def\U{{\cal U}}
\def\L{{\cal L}}

\def\mz{{m_0}}

\newtheorem{thm}{Theorem}[section]

\newtheorem{lem}[thm]{Lemma}

\newtheorem{defn}[thm]{Definition}

\newcommand{\thmref}[1]{Theorem~{\rm \ref{#1}}}
\newcommand{\lemref}[1]{Lemma~{\rm \ref{#1}}}

\newcommand{\defref}[1]{Definition~{\rm \ref{#1}}}

\newcounter{neweqn}

\newcommand{\beq}[1]{\begin{equation} \refstepcounter{neweqn} \label{#1}}

\newcommand{\eeq}{\end{equation}}
\newcommand{\bed}{\begin{displaymath}}
\newcommand{\eed}{\end{displaymath}}
\newcommand{\ben}{\begin{eqnarray*}}
\newcommand{\een}{\end{eqnarray*}}
\newcommand{\disp}{\displaystyle}

\newcommand{\bedd}{\bed\begin{array}{l}}
\newcommand{\eedd}{\end{array}\eed}
\newcommand{\al}{\alpha}

\newcommand{\sg}{\sigma}

\newcommand{\M}{{\cal M}}
\newcommand{\e}{\varepsilon}
\newcommand{\Ga}{\Gamma}

\newcommand{\wdt}{\widetilde}

\newcommand{\cd}{(\cdot)}
\newcommand{\rr}{{\mathbb R}}

\def\({\left(}
\def\){\right)}
\def\one{{\hbox{1{\kern -0.35em}1}}}

\newcommand{\bea}{\bed\begin{array}{rl}}
\newcommand{\eea}{\end{array}\eed}
\newcommand{\ad}{&\!\!\!\disp}

\newcommand{\barray}{\begin{array}{ll}}
\newcommand{\earray}{\end{array}}

\def\para#1{\vskip 0.4\baselineskip\noindent{\bf #1}}
\def\qed{$\qquad \Box$}

\newcommand{\B}{{\mathcal B}}

\title{\LARGE \bf
Existence of Saddle Points in Discrete Markov Games and Its
Application in Numerical Methods for Stochastic Differential Games}

\author{Q. S. Song and G. Yin
\thanks{Q. S. Song is with Department of Mathematics, Wayne State
University, Detroit, MI 48202, {\tt song@math.wayne.edu}.
Research of this author was supported in part
 by WSU Graduate Research Assistantship.}%
\thanks{G. Yin is with Department of Mathematics, Wayne State
University, Detroit, MI 48202, {\tt gyin@math.wayne.edu}.
Research of this author was supported in part
by the National Science Foundation.
}%
}

\begin{document}

\maketitle
\thispagestyle{empty}
\pagestyle{empty}

\begin{abstract}
This work establishes sufficient conditions for existence of
saddle points in discrete Markov games. The result
 reveals the relation between dynamic games
and  static games using dynamic programming equations. This result
enables us to prove existence of saddle points of non-separable
stochastic differential games of regime-switching diffusions under
appropriate conditions.

\end{abstract}

\section{Introduction}

The merge of differential games and regime-switching models stems
from a wide range of applications in communication networks,
complex
systems, and financial engineering. Many problems arising in, for
example, pursuit-evasion games, queueing systems in heavy traffic,
risk-sensitive control, and constrained optimization problems, can
be formulated as two-player stochastic differential games
\cite{BasarB,ElliottK,FM}. In another direction, recent applications
for better describing
the random environment leads to the use of the
so-called regime-switching models; see
\cite{Mariton,Rolski,YinK05,Zhang01,ZhouY} and
many references therein.
Since for many problems arising in applications,
closed-form solutions are difficult to obtain. As a viable
alternative, one is contended with numerical approximations
\cite{RF,SongYZ,YinSZ}. A systematic approach of numerical
approximation for
stochastic differential games was provided in \cite{Kushner02} using Markov
chain approximation methods. The major difficulty in dealing with
such game problems is to prove the existence of the value of the
game. To ensure the existence of saddle points, separability with
respect to controls for objective function and the drift of the
diffusion is required in \cite{Kushner02}. It would be nice to be
able to relax the separability condition.

Markov chain approximations of stochastic differential games are
indeed discrete Markov games. In this paper, we aim to develop
sufficient conditions for the existence of saddle point of discrete
Markov games. In the proof, we start with dynamic programming
equation together with static game results obtained by Sion
\cite{Sion} and von Neumann \cite{von}, discover the relations
between static games and dynamic games by a series of inequalities.
This approach enables us to treat non-separable discrete Markov
games with respect to controls. By virtue of results in discrete
Markov games, we can easily prove the existence of saddle points of
discrete Markov games arising in numerical approximations of
stochastic differential games when a discretization parameter $h$ is
used.
As $h\to 0$, we are able to obtain the existence of saddle points of
non-separable stochastic differential games using weak convergence
techniques in \cite{KushnerD} and \cite{Kushner02}.

The rest of the paper is arranged as follows. Section II begins with
the formulation of the discrete Markov games. Section III presents
sufficient conditions for the existence of saddle points of
discrete Markov games for both ordinary control and relaxed control
spaces, respectively. Section IV applies the results in the discrete
Markov games to stochastic differential games. Section V
concludes the paper with further remarks.

\section{Formulation}

Consider a two-player discrete Markov zero-sum game. Let $S$ be a
finite state space of a Markov chain, and $\partial S \subset S$ be
a collection of absorbing states. Control space $U_1$ and $U_2$ for
player $1$ and player $2$ are compact subsets of $\rr$. [For
notational simplicity, we have chosen to treat real-valued controls
in this paper.] Let $\{\xi_n, n<\infty\}$ be a controlled
discrete-time Markov chain, whose time-independent transition
probabilities controlled by a pair of sequences $\{(u_{1,n},
u_{2,n}), n<\infty\}$ is
\begin{equation} \label{trans1}
p(x,y|r_1,r_2) = P\{\xi_{n+1} = y| \xi_n = x, u_{1,n} =r_1, u_{2,n}
=r_2\},
\end{equation}
where $u_{i,n}\in U_i$ denote the decision at time $n$ by player
$i$.

\begin{defn} \label{adcon}
A control policy $\{(u_{1,n}, u_{2,n}), n<\infty\}$ for the chain
$\{\xi_n, n<\infty\}$ is admissible if
\begin{equation} \label{adcon-1}
P\{\xi_{n+1} = y | \xi_k, u_{1,k}, u_{2,k}, k\le n\} = p(\xi_n, y|
u_{1,n}, u_{2,n}).
\end{equation}
If there is a function $u_i\cd$ such that $u_{i,n} = u(\xi_i)$, then
we refer to $u_i\cd$ as a feedback control of player $i$.
\end{defn}

Given the running cost function $c(\cdot,\cdot,\cdot): S\times
U_1\times U_2 \mapsto \rr^+ \cup \{0\}$, and the terminal
cost function
$g\cd:S \to \rr^+ \cup \{0\}$, the cost for an initial $\xi_0 = x\in
S$ and an admissible control policy $(u_1,u_2) = \{(u_{1,n},
u_{2,n}): n<\infty\}$ is defined by
\begin{equation} \label{cost1}
W(x,u_1, u_2) = E_x^{u_1,u_2} [ \sum_{n=0}^{N-1} c(\xi_n, u_{1,n},
u_{2,n}) + g(\xi_N)],
\end{equation}
where $N = \min\{n: \xi_n\in \partial S\}$ and $E_x^{u_1,u_2}$ is
the expectation given that initial $\xi_0 = x$ and control
$(u_1,u_2)$.

In the discrete Markov game, player $1$ wants to minimize the cost,
while player $2$ wants to maximize. The two
players have different information available
depending on who makes the decision first (or who ``goes first'').
Using $\U_i(1)$ to denote the space of the admissible ordinary
controls that player $i$ goes first.
That is, for $u_i\in \U_i(1)$, there exists a sequence of
measurable functions $F_n\cd$ taking values in $U_i$
such that $u_{i,n} = F_n (\xi_k, k\le n; u_{1,k}, u_{2,k}, k<n).$
Similarly, using $\U_i(2)$ to denote the collection of the admissible
ordinary
controls that player $i$ goes last, that is,
$u_i \in \U_i(2)$ is determined by a sequence of measurable
functions $\wdt F_n\cd$ taking values in $U_i$ such that
$u_{i,n} = \wdt F_n (\xi_k, k\le n; u_{i,k}, k<n; u_{j,k}, k\le n,
j\neq i).$

To proceed,
we define upper and lower values by
\begin{equation} \label{uval}
V^+(x) = \min_{u_1\in \U_1(1)} \max_{u_2\in \U_2(2)} W(x,u_1,u_2)
\end{equation}
\begin{equation} \label{lval}
V^-(x) = \max_{u_2\in \U_2(1)} \min_{u_1\in \U_1(2)} W(x,u_1,u_2)
,\end{equation}
respectively. It is obvious $V^-(x) \le V^+(x)$ for $\forall x\in
S$.
If the lower value and upper value are equal, then we say there
exists a saddle point for the game, and its value is
\begin{equation} \label{sad}
V(x) = V^+(x) = V^-(x), \quad \forall x\in S.
\end{equation}

The corresponding dynamic programming equation is
\begin{equation} \label{udpe}
V^+(x)  = \min_{r_1\in U_1} \max_{r_2\in U_2} \{E_x[V^+(\xi_1)] +
c(x,r_1,r_2)\},
\end{equation}
\begin{equation} \label{ldpe}
V^-(x)  = \max_{r_2\in U_2} \min_{r_1\in U_1} \{E_x[V^-(\xi_1)] +
c(x,r_1,r_2)\}.
\end{equation}
Practically, we can find $V^+$ and $V^-$ in (\ref{uval}) and
(\ref{lval}) by solving (\ref{udpe}) and (\ref{ldpe}) using
iterations. This is possible owing to the following lemma. The proof
of this lemma can be found in \cite[Lemma 2]{KS69}, and a weaker
form in \cite{Zach}.
\begin{lem} \label{sdl}
$\{\xi_n, n<\infty\}$ is Markov chain with state space $S$,
absorbing states $\partial S$, and transition probability
$p(x,y|r_1,r_2)$.  Let there be a real number $\gamma>0$ with
\begin{equation} \label{sdl-1}
P(\xi_n \in \partial S| \xi_0 = x, u_{1,k}, u_{2,k}, k\le n) \ge \gamma, \quad \forall x\in S,
\end{equation}
$c(x,r_1,r_2)$ is continuous in $r_1$ and $r_2$, To each admissible control, $(u_1,u_2)$, the cost $W(x,u_1,u_2)$ is defined by (\ref{cost1}). Then $W(x,u_1,u_2)$ is finite and solutions of (\ref{udpe}) and (\ref{ldpe}) are unique. For any initial value $\{V_0^+(x): x\in S\}$, the sequence
\begin{equation} \label{iter1}
V_{n+1}^+ (x) = \min_{r_1\in U_1} \max_{r_2\in U_2} \{
E_x[V^+_n(\xi_1)] + c(x,r_1,r_2)\}
\end{equation}
converges to $V^+(x)$, the unique solution of (\ref{udpe}) as $n\to
\infty$. Analogously, for any initial $\{V_0^-(x), x\in S\}$, the
sequence
\begin{equation} \label{iter2}
V_{n+1}^- (x) = \max_{r_2\in U_2} \min_{r_1\in U_1}  \{
E_x[V^-_n(\xi_1)] + c(x,r_1,r_2)\}
\end{equation}
\end{lem}
converges to $V^-(x)$, the unique solution of (\ref{ldpe}) as $n\to
\infty$.


\section{Existence of Saddle Points}

In this section, we provide
sufficient conditions
for the existence of saddle points in discrete Markov games.
An  existence
proof is established through a series of inequalities.
In addition, the definition of relaxed controls is
given as a generalization of ordinary controls.
It is shown that saddle points always exist in relaxed control space.
\begin{defn} \label{vc}
$f(r_1,r_2)$ is said to be convex-concave with
respect to $(r_1,r_2)$, if $f(\cdot, r_2)$ is convex
and $f(r_1,\cdot)$ is concave.
\end{defn}

Next, we present a well-known minimax principle in
static games, which was obtained by Sion in \cite{Sion}.

\begin{lem} \label{minmax}
Let $M_1$ and $M_2$ be compact spaces, $\phi(\cdot,\cdot)$ be a
convex-concave function on $M_1\times M_2$, then $$\min_{r_1\in M_1}
\max_{r_2\in M_2} \phi(r_1,r_2) = \max_{r_2\in M_2} \min_{r_1\in
M_1} \phi(r_1,r_2). $$
\end{lem}

One of following two assumptions are needed for the
 existence theorem.

\begin{description}
\item [(H1)] $p(x,y|r_1,r_2)$ and $c(x,r_1,r_2)$ are
continuous and separable in $r_1$ and $r_2$.
\item [(H2)] $p(x,y|r_1,r_2)$ and $c(x,r_1,r_2)$ are
convex-concave with respect to $(r_1,r_2)$.
\end{description}

\begin{thm} \label{st} Assume
either (H1) or (H2). $\{\xi_n, n<\infty\}$ is a Markov chain as in
\lemref{sdl}.  Let $V^+(x)$ and $V^-(x)$ be associated upper and
lower values defined in (\ref{uval}) and (\ref{lval}). Then there
exists a saddle points, that is, $$V^+(x) = V^-(x), \quad \forall
x\in S.$$
\end{thm}

\para{Proof.} Define two functions $\phi^+\cd$ and $\phi^-\cd$ by
$$\phi^+(x,r_1,r_2) = \sum_{y\in S} p(x,y|r_1,r_2) V^+(y) + c(x,r_1,r_2),$$
$$\phi^-(x,r_1,r_2) = \sum_{y\in S} p(x,y|r_1,r_2) V^-(y) + c(x,r_1,r_2).$$
The dynamic programming equation of (\ref{udpe}) and (\ref{ldpe})
can be rewritten as
$$V^+(x) = \min_{r_1\in U_1} \max_{r_2\in U_2}
\{\phi^+(x,r_1,r_2)\},$$
$$V^-(x) =  \max_{r_2\in U_2} \min_{r_1\in U_1} \{\phi^-(x,r_1,r_2)\}.$$
Under either assumption (H1) or (H2), by \lemref{minmax},
\begin{equation}\label{st-1}
\min_{r_1\in U_1} \max_{r_2\in U_2} \phi^+(x,r_1,r_2) = \max_{r_2\in U_2} \min_{r_1\in U_1} \phi^+(x,r_1,r_2).
\end{equation}
Let $\rho = \max_{x\in S} \{V^+(x) - V^-(x)\} \ge 0$, then
\begin{equation} \label{st-2}
V^+(x) \le V^-(x) + \rho, \quad \forall x\in S.
\end{equation}
In particular, there exists $\hat x \in S$, so that equal holds in
(\ref{st-2}),
\begin{equation}\label{st-3}
V^+(\hat x) = V^-(\hat x) +\rho.
\end{equation}
For $\hat x$ given in (\ref{st-3}), a series of inequalities follows,
\begin{equation}\label{st-4}\barray
V^+(\hat x) \ad = \min_{r_1 \in U_1} \max_{r_2 \in U_2} \{ \phi^+(x,r_1,r_2) \} \\
\ad = \max_{r_2\in U_2} \min_{r_1\in U_1} \{ \phi^+(x,r_1,r_2)\} \\
\ad = \max_{r_2\in U_2} \min_{r_1\in U_1} \{ \sum_{y\in S} p(x,y|r_1,r_2) V^+(y) \\ \ad \hspace*{1.6in}+ c(x,r_1,r_2)\}\\
\ad \le \max_{r_2\in U_2} \min_{r_1\in U_1} \{ \sum_{y\in S} p(x,y|r_1,r_2) (V^-(y)+\rho)\\ \ad \hspace*{1.6in} + c(x,r_1,r_2)\}\\
\ad = \max_{r_2\in U_2} \min_{r_1\in U_1} \{ \phi^-(x,r_1,r_2)\} + \rho\\
\ad = V^-(\hat x) + \rho. \earray
\end{equation}
By virtue of (\ref{st-3}), we conclude all inequalities are indeed
equal in (\ref{st-4}), and this implies
$$V^+(y) = V^-(y) + \rho, \quad \forall y\in S.$$
Note that $V^+(x) = V^-(x)$ for all $x\in \partial S$. Hence
$\rho=0$. The existence of the saddle point is established. \qed

The above theorem gives sufficient conditions for the existence of
saddle points.
We note that
there always exist  saddle
points in {\it
relaxed control} space
with merely continuity assumed.

\begin{defn} \label{rc}
A control policy $\{(m_{1,n}, m_{2,n}), n<\infty\}$ for the chain
$\{\xi_n, n<\infty\}$ is said to be a relaxed control policy, if
$m_{i,n}$ is a probability measure on $\mathcal{B} (U_i)$, a
$\sg$-algebra of Borel subsets of $U_i$.
\end{defn}

More general definition of relaxed control is given by
\defref{dfn-rcon} in the context of stochastic differential games.
Let $\mathcal{P}(U_1)$ and $\mathcal{P}(U_2)$ be collection of
probability measure on $\mathcal{B} (U_1)$ and $\mathcal{B} (U_2)$.
Slightly abusing notations,  we generalize real function $f(\cdot,
\cdot)$ on $U_1\times U_2$ into a  function $f$ on ${\mathcal
P}(U_1)\times {\mathcal P}(U_2)$ as following
$$f(\mu_1, \mu_2) = \int_{U_1} \int_{U_2} f(r_1,r_2)
\mu_1(d r_1) \mu_2(d r_2).$$
Using the notation of relaxed control representation,
the transition probability function is
$$p(x,y|\mu_1, \mu_2) = \int_{U_1}\int_{U_2}
p(x,y|r_1,r_2) \mu_1(d r_1) \mu_2(d r_2), $$ and the cost
under the relaxed control policy $(m_1,m_2) = \{(m_{1,n}, m_{2,n}),
n<\infty\}$ is
$$W(x,m_1,m_2) =   E_x^{m_1,m_2} [\sum_{n=0}^{N-1} c(\xi_n, m_{1,n}, m_{2,n}) + g(\xi_N)]. $$

Using $\Gamma_i(1)$ to denote the space of admissible
relaxed controls that player $i$ goes first.
That is, for $m_i\in \Gamma_i(1)$, there exists a
sequence of measurable function $H_n\cd$ taking values
in $\mathcal{P} (U_i)$ such that
$$m_{i,n} = H_n(\xi_k, k\le n, m_{1,k}, m_{2,k}, k<n).$$
Analogously, using $\Gamma_i(2)$ to denote the space of
admissible relaxed controls that player $i$ goes last.
That is, $m_i\in \Gamma_i(2)$, there exists a sequence of
measurable function $\wdt H_n\cd$ taking values
in $\mathcal{P} (U_i)$ such that
$$m_{i,n} = \wdt H_n(\xi_k, k\le n; m_{1,k},
k<n; m_{j,k}, k\le n, j\neq i).$$
The upper and lower values associated with
relaxed control space are defined by
\begin{equation} \label{uvalr}
V^+_m(x) = \min_{m_1\in \Ga_1(1)} \max_{m_2\in \Ga_2(2)} W(x,m_1,m_2)
\end{equation}
\begin{equation} \label{lvalr}
V^-_m(x) = \max_{m_2\in \Ga_2(1)} \min_{m_1\in \Ga_1(2)}
W(x,m_1,m_2),
\end{equation}
respectively. To proceed, we present another static game result
obtained by von Neumann \cite{von}.

\begin{lem} \label{vonl}
Let $M_1$ and $M_2$ be finite sets. Let $\phi(\cdot,\cdot)$ be a
function on $M_1\times M_2$, $\mu_1\in {\cal P}(M_1)$ and $\mu_2\in
{\cal P}(M_2)$ be probability measure on $M_1$ and $M_2$, then
\begin{equation}\label{vonl-1} \barray
\ad \min_{\mu_1\in \mathcal{P}(M_1)} \max_{\mu_2\in
\mathcal{P}(M_2)} \phi(\mu_1, \mu_2) =  \\ \ad \hspace{1in}
\max_{\mu_2\in \mathcal{P}(M_2)} \min_{\mu_1\in \mathcal{P}(M_1)}
\phi(\mu_1, \mu_2). \earray \end{equation}
\end{lem}

\begin{thm} \label{rsad}
$\{\xi_n, n<\infty\}$ is a Markov chain as in \lemref{sdl} with
relaxed control used. Assume $p(x,y|\cdot,\cdot)$ and
$c(x,\cdot,\cdot)$ are continuous on $U_1\times U_2$. Let $V_m^+(x)$
and $V_m^-(x)$ be associated upper and lower values of (\ref{uvalr})
and (\ref{lvalr}). Then there always exists a saddle point, that is
$$V_m^+(x) = V_m^-(x), \quad \forall x \in S.$$
\end{thm}

\para{Proof.} Define two functions $\phi_m^+\cd$ and $\phi_m^-\cd$ by
$$\phi_m^+(x,\mu_1,\mu_2) = \sum_{y\in S} p(x,y|\mu_1,\mu_2) V_m^+(y) + c(x,\mu_1,\mu_2),$$
$$\phi_m^-(x,\mu_1,\mu_2) = \sum_{y\in S} p(x,y|\mu_1,\mu_2) V_m^-(y) + c(x,\mu_1,\mu_2).$$
Then dynamic programming equation in relaxed control space can be
written by
$$V_m^+(x) = \min_{\mu_1\in {\cal P}(U_1)} \max_{\mu_2\in
{\cal P}(U_2)} \{\phi_m^+ (x,\mu_1,\mu_2)\},$$
$$V_m^-(x) = \max_{\mu_2\in
{\cal P}(U_2)} \min_{\mu_1\in {\cal P}(U_1)} \{\phi_m^-
(x,\mu_1,\mu_2)\},$$

Note that $c(x,\cdot, \cdot)$ is continuous in compact set
$U_1\times U_2$. Hence for $\forall \e>0$, there exists a finite
subset $U_1^{\e} \times U_2^{\e} \subset U_1 \times U_2$, such that
\begin{equation}\label{rsad-2} \barray
\ad \Big |\min_{\mu_1\in \mathcal{P}(U_1)} \max_{\mu_2\in
\mathcal{P}(U_2)} c(x, \mu_1, \mu_2)   \\ \ad \hspace*{.6in} -
\min_{\mu_1^{\e}\in \mathcal{P}(U_1^{\e})} \max_{\mu_2^{\e}\in
\mathcal{P}(U_2^{\e})} c(x, \mu_1^{\e}, \mu_2^{\e})\Big | < \e.
\earray\end{equation}
\begin{equation}\label{rsad-3} \barray
\ad \Big |\max_{\mu_2\in \mathcal{P}(U_2)}
\min_{\mu_1\in \mathcal{P}(U_1)} c(x, \mu_1, \mu_2)   \\
\ad \hspace*{.6in} - \max_{\mu_2^{\e}\in \mathcal{P}(U_2^{\e})}
\min_{\mu_1^{\e}\in \mathcal{P}(U_1^{\e})} c(x, \mu_1^{\e},
\mu_2^{\e})\Big | < \e. \earray\end{equation}

Forcing to the limit as $\e\to 0$ in (\ref{rsad-2}) and
(\ref{rsad-3}), as well as using \lemref{vonl}, we have
\begin{equation}\label{rsad-4} \barray
\ad \min_{\mu_1\in \mathcal{P}(U_1)} \max_{\mu_2\in
\mathcal{P}(U_2)} c(x, \mu_1, \mu_2) =  \\ \ad
\hspace*{.8in}\max_{\mu_2\in \mathcal{P}(U_2)} \min_{\mu_1\in
\mathcal{P}(U_1)} c(x, \mu_1, \mu_2). \earray\end{equation}

Similarly, we obtain equality for function $p(x,y|\cdot,\cdot)$,
\begin{equation}\label{rsad-5} \barray
\ad \min_{\mu_1\in \mathcal{P}(U_1)} \max_{\mu_2\in
\mathcal{P}(U_2)} p(x,y| \mu_1, \mu_2) =  \\ \ad
\hspace*{.8in}\max_{\mu_2\in \mathcal{P}(U_2)} \min_{\mu_1\in
\mathcal{P}(U_1)} p(x,y| \mu_1, \mu_2). \earray\end{equation}

Equalities in (\ref{rsad-4}) and (\ref{rsad-5}) implies
\begin{equation}\label{rsad-1} \barray
\ad \min_{\mu_1\in \mathcal{P}(U_1)}  \max_{\mu_2\in
\mathcal{P}(U_2)} \phi_m^+(x, \mu_1, \mu_2) =  \\ \ad
\hspace*{.8in}\max_{\mu_2\in \mathcal{P}(U_2)} \min_{\mu_1\in
\mathcal{P}(U_1)} \phi_m^+(x, \mu_1, \mu_2). \earray\end{equation}
The rest of this proof is similar to the lines of inequalities
(\ref{st-4}). The details are omitted. \qed

\section{Numerical Methods Regime-Switching Stochastic Differential Games}
In this section, we formulate stochastic differential games with
regime switching. Numerical methods using Markov chain approximation
leads to a sequence of discrete Markov games discussed in the
previous section. The use of \thmref{st} gives sufficient conditions
for the existence of saddle points,
and facilitates the
proof.

\subsection{Formulation}
Consider a two-player stochastic game of regime-switching
diffusions. For a finite set $\M=\{1,\ldots,\mz\}$, $x\in \rr^\lz$,
$b(\cdot,\cdot,\cdot): \rr^\lz \times \M \times \rr \times
\rr\mapsto
 \rr^\lz$, $\sg(\cdot,\cdot,\cdot):
  \rr^\lz \times \M \mapsto \rr^\lz \times
\rr^\lz$, the dynamic system is given by
\begin{equation} \label{mod1} \barray
\ad x(t) = x(0) + \int_0^t b(x(s), \al(s), u_1(s),u_2(s)) ds \\
\ad \hspace{1.3in} +  \int_0^t \sg(x(s),\al(s)) d w(s), \earray
\end{equation}
where for each $i=1,2$, $u_i(\cdot)$ is a control for player $i$,
$w(\cdot)$ is a standard $\rr^\lz$-valued Brownian motion, and
$\al\cd$ is a continuous-time Markov chain having state space $\M$
with generator $Q=(q_{\iota,\ell}) \in \rr^{\mz\times \mz}$. Let
$\{\F_t: 0\le t\}$ be a filtration, which might depend on controls,
and which measures at least $\{(w(s),\al(s)): s \le t\}$. We suppose
that for each $i=1,2$, $u_i(\cdot)$ is $\F_t$-adapted taking values
in a compact subset $U_i \subset\rr$, which are called
\textit{admissible controls}. Denote  $A(x,\iota) = \sg(x,\iota)
\sg'(x,\iota)=(a_{j_0k_0}(x,\iota)) \in \rr^\lz \times \rr^\lz$,
which is symmetric and positive definite.

Let $G \subset \rr^\lz$ be a compact set that is the closure of its
interior $G^0$ and $\tau$ be the first exit time of $x(t)$ from
$G^o$ with
\begin{equation} \label{tau-def}
\tau = \min\{t:x(t)\notin G^o\}.
\end{equation}
Using a real number $\beta > 0$ to denote the discount factor, let
the cost function  be
\begin{equation} \label{costfun1} \barray
 W(x,\ad  \iota,u) = E_{x,\iota}^{u}\Big[ \int_0^\tau
e^{-\beta s} \wdt{k}(x(s),\al(s),u(s)) ds \\ \ad \hspace{1.5in}+
\wdt{g}(x(\tau),\al(\tau))\Big],\earray
\end{equation}
where $\wdt{k} \cd$ and $\wdt g\cd$ are functions representing the
running cost and terminal cost, respectively, and
 $E^u_{x,\iota}$ denotes the expectation taken with
the initial data $x(0)=x$ and $\al(0)=\iota$ and given control
process $u(\cdot) = (u_1\cd, u_2\cd)$. Next, we introduce the
relaxed control representation; see \cite{Kushner02, KushnerD}.

\begin{defn}\label{dfn-rcon}
{\rm Let ${\cal B}({U} \times[0,\infty))$ be the $\sigma$-algebra of
Borel subsets of ${U} \times [0,\infty)$. An {\it admissible relaxed
control}  $m(\cdot)$ is a measure on ${\cal B}({U}
\times[0,\infty))$ such that $m({U}\times[0,t]) = t$ for each $t \ge
0$. Given a relaxed control $m(\cdot)$, there is an $m_t(\cdot)$
such that $m(drdt) = m_t(dr)dt$. In fact, we can define $m_t(B) =
\lim_{\delta \to 0} \disp{\frac {m(B\times[t-\delta,t])} {\delta}}$
for $B\in \B(U)$.}\end{defn}

To proceed, we need the following assumptions.

\begin{itemize}
\item [(A1)] For each $\iota \in \M$,
$\wdt{k} (\cdot, \iota, \cdot, \cdot)$ and
 $b(\cdot, \iota, \cdot, \cdot)$  are continuous functions
 on the compact set $G \times U_1 \times U_2$.

 \item [(A2)] For each $\iota\in \M$, the functions $\sg (\cdot, \iota)$
 and $\wdt{g} (\cdot, \iota)$ are continuous on $G$.

\item [(A3)] Equation (\ref{mod1}), where the
controls are replaced by relaxed
controls,
has a unique weak sense
 solution (i.e., unique in the sense of in distribution)
 for each admissible triple $(w\cd, \al\cd, m\cd)$, where
 $m\cd = (m_1\cd, m_2\cd)$.
 \item [(A4)] For any $\iota\in \M,$ $j_0,k_0 \in \{1,2,\ldots,\lz\},
 j_0\neq k_0,$ $a_{j_0j_0}(x,\iota) > \sum_{k_0\neq j_0}
 |a_{j_0k_0} (x,\iota)|$.
 \item [(A5)] Let
 $\hat{\tau}(\phi)  =
 \left\{  \barray
\ad \infty, \ \hbox{ if  } \ \phi(t) \in
 G^o \ \hbox{ for all  } \ t<\infty, \\
\ad \inf \{ t: \phi(t) \notin G^o\} \
 \hbox{ otherwise. } \earray
 \right.$
 The function $\hat \tau \cd$
 is continuous as a mapping from $D[0,\infty)$ to $[0,\infty]$
 with probability one relative to the measure
 induced by any solution with initial condition
 $(x,\iota)$,
 where $D[0,\infty)$ denotes the space of
 functions that are right continuous and have left limits endowed
 with the Skorohod topology, and $[0,\infty]$ is the
 interval $[0,\infty)$
 compactified (see \cite[p. 259]{KushnerD}).

\item [(A6)] The functions $b\cd$ and $\wdt k\cd$ are separable
 in $r_1$ and $r_2$ for every $(x,\iota)\in G\times \M$. That is,
 $b(x,\iota,r_1,r_2) = \sum_{i=1}^2 b^i(x,\iota,r_i)$ and
 $\wdt k (x,\iota,r_1,r_2) = \sum_{i=1}^2 \wdt k^i(x,\iota,r_i)$.

\item[(A7)] The cost $\wdt k \cd$ is convex-concave with respect to
$(r_1,r_2)$, and there exist $\rr^\lz$-valued continuous functions
$b^i(x,\iota)$ ($i=0,1,2,3$) such that
$b(x,\iota, r_1, r_2)= r_1 r_2 b^0(x,\iota) + r_1 b^1(x,\iota)
+ r_2 b^2(x,\iota)+ b^3(x,\iota).$
\end{itemize}

Assumption (A4) is used for construction of transition probabilities
of the approximating Markov chain.  It requires that the diffusion
matrix be diagonally dominated. If the given dynamic system does not
satisfy (A4), then we can adjust the coordinate system to satisfy
assumption (A4); see \cite[p. 110]{KushnerD}. (A5) is a broad
condition that is satisfied in most applications. The main purpose
is to avoid the {\it tangency} problem discussed in \cite[p.
278]{KushnerD}. Later, we will establish the existence of {\it
saddle points} using either (A6) or (A7) in addition to (A1)--(A5).
Condition (A7) allows non-separable differential games with respect
to controls.

Now we are ready to define upper values, lower values, and
\textit{saddle points} of differential games; see \cite{Kushner02}
for the corresponding definitions of systems without regime
switching. Let $\U_i$ be collection of all admissible ordinary
control with respect to $(w\cd, \al\cd)\}$. For $\Delta>0$, Let
$\U_i(\Delta) \subset \U_i$ such that $u_i\cd$ are piecewise
constant on the intervals $[k\Delta, k\Delta+\Delta), k =
0,1,2,\ldots$, and $u_i(k\Delta)$ is $\F_{k\Delta}$-measurable.

Let $\L_1(\Delta)\subset \U_1(\Delta)$ denote the set of such
piecewise constant controls for player $1$ that are determined by
measurable real-valued functions $Q_{1,n}(\cdot)$
\begin{equation} \label{defcon1}
u_1(n\Delta)  = Q_{1,n}(w(s),\al(s), u(s), s<n\Delta),
\end{equation}
We can define $\L_2(\Delta)$ and the associated rule $u_2$ for
player $2$ analogous to (\ref{defcon1}).

Thus we can always suppose that if the control of (for example)
player $1$ is determined by a form such as (\ref{defcon1}). Then (in
relaxed control terminology) the law of $(w(t), \al(t), m_2(t))$ for
$n\Delta \le t < (n+1)\Delta$ is determined recursively by past
information
 \begin{equation} \label{defcon3}
 \{w(s), \al(s), m_2(s), s<t, , m_1(s), s \le n\Delta \}.
 \end{equation}

\begin{defn}\label{defval1}{\rm
 For initial condition $x(0) = x, \al(0) = \iota$, define the upper
 and lower values for the game as
 \begin{equation} \label{defupval1}
 V^+(x,\iota) = \lim_{\Delta\to 0} \inf_{u_1 \in \L_1(\Delta)}
 \sup_{u_2\in \U_2} W(x, \iota, u_1, u_2),
 \end{equation}
 \begin{equation} \label{deflwval1}
 V^-(x,\iota) = \lim_{\Delta\to 0} \sup_{u_2 \in \L_2(\Delta)}
 \inf_{u_1\in \U_1} W(x, \iota, u_1, u_2).
 \end{equation}
 If the lower and upper value are equal, then we say there
 exists a saddle point for the game, and its value is
 \begin{equation} \label{sadpt1}
 V^+(x,\iota) = V^-(x,\iota) = V(x,\iota),\, \forall x\in G, \iota
 \in \M.
\end{equation}
}
\end{defn}

\subsection{Markov Chain Approximations}
Here, we will construct a two-component Markov chain. The
discretization of differential game leads to a sequence of discrete
Markov games. The approximation is of finite difference type. The
basis of the approximation is a discrete-time, finite-state,
controlled Markov chain $\{(\xi_n^h, \al_n^h): n<\infty \}$ whose
properties are {\it locally consistent} with that of (\ref{mod1}).

For each $h>0$, let $G_h$ be a finite subset of $G$ such that
$d(G_h,G)\to 0$ as $h\to 0$, where $d\cd$ is a metric defined by
\begin{equation} \label{met} d(G,G_h) = \max_{p\in G} \min_{q\in G_h} d(p,q).
\end{equation} Let $\{(\xi_n^h, \al_n^h): n<\infty \}$ be a controlled
discrete-time Markov chain on a discrete state space $G_h \times \M$
with transition probabilities denoted by
$p^h((x,\iota),(y,\ell)|r)$, where $r=(r_1,r_2)\in U_1\times U_2$.
We use $(u_{1,n}^h, u_{2,n}^h)$ to denote the actual control action
for the chain at discrete time $n$. Suppose we have a positive
function $\Delta t^h\cd$ on $G_h\times \M \times U_1 \times U_2$
such that $\sup_{x,\iota,r} \Delta t^h(x,\iota,r) \to 0$ as $h\to
0$, but $\inf_{x,\iota,r} \Delta t^h(x,\iota,r) > 0$ for each $h>0$.
We take an interpolation of the discrete Markov chain $\{(\xi_n^h,
\al_n^h)\}$ by using interpolation interval $\Delta t_n^h = \Delta
t_n^h(\xi_n^h, \al_n^h, u_{1,n}^h, u_{2,n}^h)$. Now we give the
definition of local consistency.

\begin{defn}\label{deflc1}
{\rm Let $\{p^h((x,\iota),(y,\ell)|r)\}$ for $(x,\iota)$ and $(y,
\ell)$ in $G_h\times \M$ and $r\in U_1 \times U_2$ be a collection
of well-defined transition probabilities for the two-component
Markov chain $\{(\xi_n^h,\al_n^h)\}$, approximation to
$(x(\cdot),\al(\cdot))$. Define the difference $\Delta \xi_n^h =
\xi_{n+1}^h -\xi_n^h$. Assume $\lim_{h\to 0} \sup_{x,\iota,r} \Delta
t^h (x,\iota,r) = 0$. Denote by $E_{x,\iota,n}^{r,h}$, ${\rm
cov}_{x,\iota,n}^{r,h}$ and $p_{x,\iota,n}^{r,h}$  the conditional
expectation, covariance, and probability given $\{\xi_k^h, \al_k^h,
u_{1,k}^h, u_{2,k}^h, k\le n, \xi_n^h = x, \al_n^h = \iota,
(u_{1,n}^h, u_{2,n}^h) = r\}$. The sequence $\{(\xi_n^h,\al_n^h)\}$
is said to be \textit{locally consistent} with (\ref{mod1}), for
$\Delta t^h = \Delta t^h(x,\iota,r)$, if
\begin{equation} \label{deflc2} \barray
\ad E_{x,\iota,n}^{r,h} \Delta \xi_n^h = b(x,\iota,r) \Delta t^h +
o(\Delta t^h), \\ \ad {\rm cov}_{x,\iota,n}^{r,h} \Delta \xi_n^h =
A(x,\iota) \Delta t^h + o(\Delta t^h),
\\
\ad p_{x,\iota,n}^{r,h}\{\al_{n+1}^h = \ell\} = q_{\iota\ell} \Delta
t^h + o(\Delta t^h),  \textrm{ for } \ell\neq \iota,\\
\ad p_{x,\iota,n}^{r,h}\{\al_{n+1}^h = \iota\} =  (1 +
q_{\iota\iota}) \Delta t^h + o(\Delta t^h), \\
\ad \sup_{n, \omega \in \Omega} |\Delta \xi_{n}^h| \to 0 \ \hbox{ as
} \ h\to 0, \earray
\end{equation}
}
\end{defn}

To approximate the cost defined in (\ref{costfun1}), we define a
cost function using the Markov chain above. Let $$t_n^h =
\sum_{j=0}^{n-1} \Delta t_j^h \hbox{ and } N_h = \inf\{n: \xi_n^h
\notin G_h^o\}.$$ The cost for $u^h = \{(u_{1,n}^h, u_{2,n}^h)\}$
and initial $(x,\iota)$ is
\begin{equation} \label{costfun2} \barray
\ad W^h(x,\iota, u^h) = E_{x,\iota} \Big[\sum_{n=0}^{N_h-1}
e^{-\beta t_n^h} \Delta t_n^h \cdot \\ \ad \hspace{.5in}
\wdt{k}(\xi_n^h, \al_n^h, u_{1,n}^h, u_{2,n}^h)
 + \wdt{g}(\xi_{N_h}^h,
\al_{N_h}^h)\Big], \earray
\end{equation}

Using $\U_i^h(1)$ to denote the space of the ordinary controls that
player $i$ goes first, and its strategy is defined by measurable
functions of the type similar to (\ref{defcon1}). That is, for
$u^h_i \in \U_i^h(1)$, $u_{i,n}^h$ is determined by
$$\{ \xi_k^h, \al_k^h, k\le n; u_{1,k}^h, u_{2,k}^h, k<n\}.$$
By $\U_i^h(2)$ denote the collection of the ordinary controls that
player $i$ goes last. For $u_i^h\in \U_i^h(2)$, $u_{i,n}^h$ is
determined by
$$\{\xi_k^h, \al_k^h, k\le n; u_{i,k}^h, k<n; u_{j,k}^h, k\le n,
j\neq i\}.$$
The associated upper and lower values  is defined as
\begin{equation} \label{defupval2}
V^{h,+} (x,\iota) = \inf_{u^h_1\in \U^h_1(1)} \sup_{u^h_2\in
\U^h_2(2)} W^h(x,\iota,u^h_1,u^h_2),
\end{equation}
\begin{equation} \label{deflwval2}
V^{h,-} (x,\iota) = \sup_{u^h_2\in \U^h_2(1)} \inf_{u^h_1\in
\U^h_1(2)} W^h(x,\iota,u^h_1,u^h_2).
\end{equation}

\subsection{
Saddle Points for the Markov Chain
Approximation}
In this section, we present a local consistent discrete Markov game
of $\{(\xi_n^h, \al_n^h)\}$ generated by central finite difference
scheme for analysis purpose. Under assumptions (A1)--(A5) together
with either (A6) or (A7), we can apply \thmref{st} to show the
existence of saddle points for each $h$. By forcing the limit $h\to
0$, the upper (lower) values converge to that of stochastic
differential game by \lemref{lem:conv}, and it results in the
existence of saddle points.

First, the transition probabilities for $\{(\xi_n^h, \al_n^h)\}$ are
\begin{equation} \label{transpb3}\barray
\ad p^h((x,\iota),(x\pm e_{j_0}h,\iota)|r) = \\ \ad \quad \frac {\pm
h b_{j_0}(x,\iota,r) +  a_{{j_0}{j_0}}(x,\iota)- \sum_{{k_0}\neq
{j_0}}|a_{{j_0}{k_0}}(x,\iota)|} {2(D^h(x,\iota) - \beta h^2)},\\
\ad \hspace{1.3in} \textrm{for }
{j_0}=1,2,\ldots,\lz,\\
\ad p^h((x,\iota),(x+e_{j_0} h + e_{k_0} h ,\iota)|r) = \frac {1/2
\cdot a_{{j_0}{k_0}}^+(x,\iota)} {D^h(x,\iota) - \beta h^2}, \\ \ad
p^h((x,\iota),(x - e_{j_0} h - e_{k_0} h ,\iota)|r) =  \frac {1/2
\cdot a_{{j_0}{k_0}}^+(x,\iota)} {D^h(x,\iota) - \beta
h^2}, \\
\ad  \hspace{1.5in} \textrm{for }  {j_0}<
k_0\\
\ad p^h((x,\iota),(x + e_{j_0} h - e_{k_0} h ,\iota)|r) = \frac
{1/2\cdot a_{j_0k_0}^-(x,\iota)} {D^h(x,\iota) - \beta h^2},\\
\ad  \hspace{1.5in} \textrm{for }  j_0\neq k_0,\\
\ad p^h((x,\iota),(x, \ell)|r) = \frac{q_{\iota\ell}
h^2}{D^h(x,\iota) - \beta h^2},  \ \ell\neq
\iota, \\
\ad p^h((x,\iota),(y, \ell)|r) = 0, \hbox{ otherwise.}\earray
\end{equation}
where $$D^h(x,\iota ) = \sum_{{j_0}=1}^\lz a_{{j_0}{j_0}}(x,\iota) -
\sum_{{j_0}< {k_0}} |a_{{j_0}{k_0}}(x,\iota)| -q_{\iota\iota} h^2 +
\beta h^2,$$ Set the interpolation interval as
$\Delta t^h(x,\iota) = h^2/D^h(x,\iota). $
By (A4), $D^h(x,\iota) - \beta h^2 > 0$. Also, we have $\sum_{(y,
\ell)  } p^h((x,\iota),(y,\ell) | r) = 1$. To ensure that $p^h\cd$
is always nonnegative, we require
\begin{equation} \label{asmp-h} h\le \frac {\min_{j_0} \{a_{j_0j_0}(x,\iota) -
\sum_{k_0\neq j_0} |a_{j_0 k_0}(x,\iota)|\}} {\max_r
|b_{j_0}(x,\iota,r_1,r_2)|}.\end{equation}

\begin{lem}\label{cst3}
Assume {\rm (A1), (A2)}, {\rm(A4)}, and $h$ satisfies
{\rm(\ref{asmp-h})}. The Markov chain $(\xi^h_n,\al^h_n)$ with
transition probabilities $\{p^h(\cdot)\}$ and interpolation
 $\Dl t^h\cd$ defined above is locally consistent with
{\rm(\ref{mod1})}.
\end{lem}
\para{Proof.}
The criterion in (\ref{deflc2}) can be verified through a series of
calculations, thus details are omitted. \qed

\begin{thm} \label{sadptthm1}
Assume {\rm(A1)--(A5)}, either  {\rm(A6)} or {\rm(A7)}, and $G_h$ is
a finite set defined above (\ref{met}). For $x\in G_h$ and $\iota
\in \M$, a Markov chain is defined by {\rm(\ref{transpb3})} .  Let
$V^{h,+}(x,\iota)$ and $V^{h,-}(x,\iota)$ be the associated upper
and lower values defined in
{\rm(\ref{defupval2})} and {\rm(\ref{deflwval2})} in the control
spaces $\U_i^h(1)$ and $\U_j^h(2)$ . Then there  exists a saddle
point
\begin{equation} \label{sadpt2}
 V^{h,+}(x,\iota) = V^{h,-}(x,\iota),
\end{equation}
provided $h$ satisfies {\rm(\ref{asmp-h})}.
\end{thm}
\para{Proof.} The contraction condition (\ref{sdl-1}) satisfies for
the discount factor $\beta>0$. Let
$$p((x,\iota)(y,\ell)|r_1,r_2) = p^h(((x,\iota)(y,\ell)|r_1,r_2),$$
$$c(x,\iota, r_1,r_2) = e^{-\beta \Delta t^h(x,\iota)} \Delta
t^h(x,\iota) \wdt k (x,\iota, r_1,r_2).$$ Assumptions (A6) and (A7)
lead to (H1) and (H2), respectively. The result holds applying
\thmref{st}.\qed

Although the proof of next lemma is rather complicated and not
trivial, the proof is referred to weak convergence techniques in
\cite{KushnerD}, \cite{Kushner90}, and \cite{Kushner02} due to the
limit of space.
\begin{lem}  \label{lem:conv}
Assume that the conditions of \thmref{sadptthm1}  are satisfied.
Then for the approximating Markov chain, we have
 \begin{equation} \label{sadpt6}
  \lim_{h\to 0}V^{h,+}(x,\iota)= V^+(x,\iota),
 \end{equation}
  \begin{equation} \label{sadpt5}
  \lim_{h\to 0}V^{h,-}(x,\iota)= V^-(x,\iota).
 \end{equation}

\end{lem}

\begin{thm} \label{sadptthm4}
Assume the conditions of \thmref{sadptthm1}  are satisfied.
 Then
the differential game has saddle point in the sense
\begin{equation} \label{sadpt3}
 V^+(x,\iota) = V^-(x,\iota).
\end{equation}
\end{thm}

\section{Further Remarks}

The key part of zero-sum game problems is existence of saddle point.
This paper is devoted to sufficient condition for the existence of
saddle point in discrete Markov game. Using dynamic programming
equation method, we are able to use static game results of Sion
\cite{Sion} and von Neumann \cite{von} to discover the sufficient
conditions. A direct application is  numerical methods for
stochastic differential game problems.

The transition probabilities used in (\ref{transpb3}) requires
restriction (\ref{asmp-h}) on $h$. Practically, we develop the
transition probabilities by upward finite difference scheme, so that
the generated one is well defined without restriction on $h$. It can
be routinely calculated to verify the local consistency. This kind
of discrete Markov game might have different upper and lower values
for some $h$. However, both the upper and lower values in this
situation converge to the original saddle point of differential game
$V(x)$ by \lemref{lem:conv} and \thmref{sadptthm4}. Numerical
examples in pursuit-evasion games are omitted due to the space
limit, although the numerical results clearly verify our works.

For a regime-switching system in which the Markov chain has a large
state space, we may use the ideas of two-time-scale approach
presented in \cite{YinZ}  (see also \cite{YinZ05} and references
therein) to first reduce the complexity of the underlying system and
then construct numerical solutions for the limit systems. Optimal
strategies of the limit systems can be used for constructing
strategies of the original systems leading to near
optimality.


\end{document}